\documentclass[12pt]{article}

\usepackage[margin=1in]{geometry}
\usepackage{setspace}
\usepackage{fontspec}
\usepackage{amsmath,amssymb,amsthm}
\usepackage{xurl}
\usepackage{hyperref}
\usepackage[
  backend=biber,
  style=numeric,
  sorting=none,
  maxbibnames=99,
  giveninits=true,
  doi=true,
  url=true,
  isbn=false
]{biblatex}
\usepackage[nameinlink,capitalise]{cleveref}
\usepackage{xcolor}
\usepackage{booktabs}
\usepackage{array}
\usepackage{enumitem}
\usepackage[nopatch=footnote]{microtype}
\usepackage{fvextra}

\newcommand{\code}[1]{\texttt{\small\detokenize{#1}}}
\fvset{
  fontfamily=tt,
  fontsize=\small,
  breaklines=true,
  breakanywhere=true,
  framesep=2mm,
}

\hypersetup{
  colorlinks=true,
  linkcolor=blue!60!black,
  citecolor=blue!60!black,
  urlcolor=blue!60!black,
  breaklinks=true,
  pdfborder={0 0 0},
}

\newcommand{\eps}{\varepsilon}

\newcommand{\dcode}{d_C}
\newcommand{\mc}[1]{\mathcal{#1}}

\newcommand{\bspan}[1]{\langle #1 \rangle}

\newcommand{\lean}[1]{\code{#1}}
\newcommand{\leanfile}[1]{\nolinkurl{#1}}
\newcommand{\leanat}[3]{%
  \ifnum#2=#3\relax
    l.~#2%
  \else
    ll.~#2--#3%
  \fi
}

\title{A machine-checked proof of the Dong--Yang classification of
optimal $(n,4)$ binary codes for BSCs}

\author{%
  Shenghao Yang\textsuperscript{*,\dag} \and
  Yanyan Dong\textsuperscript{*,\dag}%
}
\date{}

\begin{document}
\maketitle
\begingroup
\renewcommand{\thefootnote}{\fnsymbol{footnote}}
\footnotetext[1]{School of Science and Engineering, The Chinese University of Hong Kong, Shenzhen, Guangdong 518172, P.~R.~China.}
\footnotetext[2]{This work was supported in part by
the National Natural Science Foundation of China (NSFC) under Grants 62501515,
62171399 and 12141108, and by the Shenzhen Fundamental Research Program under Grant JCYJ20241202124023031.}
\endgroup

\begin{abstract}
We present a machine-checked Lean~4 formalization of Dong and Yang's
classification of optimal finite-length $(n,4)$ binary block codes for
binary symmetric channels.
The formalization was developed mainly by feeding the paper's proofs to an AI
tool.
To establish correctness, the authors verified the main theorem statements in
Lean and the accepted axioms.
This note discusses the corrections and simplifications made to the AI-generated formalization, and records discrepancies found in the paper during the formalization.
The Lean code is available at \url{https://github.com/shhyang/n4code_lean}.
\end{abstract}

\tableofcontents
\bigskip

\section{Introduction}
\label{sec:intro}

Many mathematical proofs are long and complex, and it is difficult to verify their correctness by human.
Machine-checked formalization promises correctness and reproducibility,
but until a few years ago even medium-scale developments typically required
years of effort by a dedicated team: proof size, specialized expertise, and
limited tooling all acted as bottlenecks.
As examples, Flyspeck~\cite{flyspeck} for the Kepler conjecture and the
formal proof of the Feit--Thompson odd-order theorem~\cite{FT} were each a
multi-year, multi-person effort combining hand-written proof structure with
extensive machine-generated case analysis. This landscape has shifted further with large language models (LLMs) that autoformalize human-written mathematics~\cite{autoformalize,formalizeassistant}, a task now within reach of coding agents such as OpenAI Codex.

Information and coding theory already has a substantial formalization
literature in Coq, the \texttt{infotheo} library~\cite{infotheo}, and 
Lean~4, the \texttt{InformationTheory} namespace in Mathlib~\cite{mathlib}.
There are also standalone Lean~4 projects formalizing $q$-ary covering
codes~\cite{coveringcodes}, self-dual constructions~\cite{selfdual}, and
Reed--Solomon and distance machinery for succinct proofs and agreement
tests~\cite{succinct,splat,zklinalg,sumcheck}. 
Lean-QIT~\cite{leanqit} provides reusable infrastructure for quantum Shannon
theory, including source and channel coding theorems.
These efforts focus on foundational definitions and classical or structural
results rather than end-to-end formalizations of contemporary research papers.

\subsection{Formalization of the Dong--Yang Classification}

We report on a machine-checked formalization in Lean~4~\cite{Lean4} of Dong and Yang's classification of optimal
finite-length $(n,4)$ binary block codes for binary symmetric
channels~\cite{DY25}.
Its main results are classification theorems (Theorems~1--5).
Theorem~1 states that for every $n$ and $\eps$ there is an optimal $(n,4)$
code that is linear or Class-I; Theorem~3 states that for $n>3$ every
optimal code is equivalent to a linear, Class-I, or Class-II code; and
Theorem~5 characterizes all optimal codes for $2\le n\le 8$.
The proofs rely on comparison techniques built around ``column
changes''---one- and two-column flips and the machinery that certifies
improvements---and span about 27 pages in the \emph{IEEE Trans.\ Inf.\ Theory}
double-column format (Table~\ref{tab:proof-hierarchy}).

We fed the paper's \LaTeX\ source to OpenAI Codex (models
\texttt{deepseek-v4-flash}, later \texttt{deepseek-v4-pro}), which first
produced definitions and a catalogue of statements to prove, each initially
carrying a \texttt{sorry}, Lean's placeholder for a missing proof.
AI agents then generated proofs to discharge the \texttt{sorry}s, and the project compiled without \texttt{sorry} in about two days.
The initial AI-assisted formalization required only light human intervention. We mainly reviewed the development plan and changed the LLM model to \texttt{deepseek-v4-pro} when the \texttt{deepseek-v4-flash} model failed to complete the formalization.

Although the Lean kernel checks every proof script, a \texttt{sorry}-free
build is not, by itself, a complete correctness argument: one must also
confirm that the Lean statements match the paper's theorems and that the
proofs depend only on trusted axioms~\cite{Lean4,flyspeck}.
This report complements the formalization with such an audit.
\Cref{sec:encoding} introduces the code model and the Lean formulations of
Theorems~1--5.
For Theorems~3--5, we verified the direct paper counterparts and additionally
proved stronger domination statements.

\begin{table}[t]
  \centering
  \caption{Hierarchical dependence among the paper's main results (Theorems~1--5).}
  \label{tab:proof-hierarchy}
  \small
  \begin{tabular}{@{}llp{0.42\textwidth}@{}}
  \toprule
  \textbf{Results} & \textbf{Technique} & \textbf{Ingredients} \\
  \midrule
  Theorem~1 & inclusive covering & Theorems~6, 8, 11; Corollaries~7, 9, 10 \\
  Theorem~2 & linear-code classification & Theorem~12; Corollary~13 \\
  Theorem~3 & all optimal codes ($n>3$) & Theorems~1, 2, 6, 8, 11; Lemma~14 \\
  Theorem~4 & linear under Class-I hypothesis & Theorems~3, 16; Lemma~14 \\
  Theorem~5 & finite classification ($2\le n\le 8$) & Theorems~2, 4, 6, 8, 11, 16, 17; Lemma~15; $n=2,3$ base cases \\
  Theorems~6, 8 & one-column analysis & Lemmas~18, 19; Theorem~20; Corollary~21 \\
  Theorem~11 & two-column analysis & Lemmas~22--27 \\
  Theorem~12 & linear-code analysis & --- \\
  Theorems~16, 17 & Class-I analysis & --- \\
  \bottomrule
  \end{tabular}
\end{table}

\Cref{sec:ai-assisted-formalization} recounts the AI-assisted development and
audits the trusted axioms, including \texttt{native\_decide} usage.
That audit also revealed redundant proofs in the generated code. We then interact with the AI tools to simplify the proofs, which is recorded in \Cref{sec:simplifications}.
The final library comprises fourteen Lean modules totaling 37{,}561 lines
(\Cref{tab:modules}) and is available at
\url{https://github.com/shhyang/n4code_lean}.

In addition to correctness, formalization helped clarify the paper's
proofs, which are written in natural language and are not always precise,
explicit, or correct as printed.
During formalization, the AI assistant often exposed these gaps: it corrected
errors, sharpened imprecise statements, and made implicit steps explicit.
We catalogue the discrepancies in \cref{sec:discrepancies}.

Within the Lean~4 ecosystem, mathlib~\cite{mathlib} already provides Hamming
distance and linear codes over finite fields
(\texttt{Mathlib/CodingTheory/}).
Although we rely on mathlib for general infrastructure, we do not build on
its coding-theory library, which targets the construction and analysis of
linear codes via generator matrices, subspaces, and minimum distance.
The Dong--Yang proofs instead use column types, equivalence under column
changes, and nonlinear codes, so we built a complementary framework for
Hamming distance, distance distributions, and maximum-likelihood decoding
on $(n,4)$ codes.

\begin{table}[t]
  \centering
  \caption{Modules, sizes, and content as of the checked and improved
  formalization.}
  \label{tab:modules}
  \small
  \begin{tabular}{@{}lrp{0.62\textwidth}@{}}
  \toprule
  \textbf{Module} & \textbf{Lines} & \textbf{Content} \\
  \midrule
  \leanfile{Definitions.lean} & 449 &
    Codes as columns; Hamming/\(\lambda\); classes; equivalence \\
  \leanfile{Basic.lean} & 610 &
    Weight/distance lemmas; count identities; \(\alpha\) basics \\
  \leanfile{Performance.lean} & 449 &
    Monotonicity of \(\lambda\); invariance under equivalence \\
  \leanfile{Comparison.lean} & 101 &
    Generic comparison engine (Lemma~18: \texttt{compare\_bij}) \\
  \leanfile{Compare.lean} & 1479 &
    One-column engine (Lemmas~19--20, Theorem~20, Corollary~21) \\
  \leanfile{ZeroColumn.lean} & 2933 &
    Replacing type-0 columns (Theorem~6, Corollary~7) \\
  \leanfile{TwoColumn.lean} & 3361 &
    Two-column \(Z\)-regions and Theorem~11 \\
  \leanfile{Reduction.lean} & 5696 &
    One-bit flip and reductions (Theorems~1,~8; Corollaries~9--10; Lemmas~14--15) \\
  \leanfile{Linear.lean} & 10775 &
    Linear comparison and classification (Theorems~2,~12; Corollary~13; Lemma~28) \\
  \leanfile{ClassI.lean} & 6930 &
    Class-I analysis (Theorems~16--17) \\
  \leanfile{FiniteN.lean} & 2352 &
    Finite classification for \(2\le n\le 8\) (Theorem~5) \\
  \leanfile{Nbig.lean} & 2246 &
    Covering for \(n>3\) (Theorems~3--4) \\
  \leanfile{Statements.lean} & 100 &
    Paper-label catalogue of theorem statements \\
  \leanfile{AxiomCheck.lean} & 80 &
    \texttt{\#print axioms} audit of headline theorems \\
  \midrule
  Total & 37{,}561 & \\
  \bottomrule
  \end{tabular}
  \end{table}

\subsection{Conventions}
We adopt the following conventions unless stated otherwise:
\begin{itemize}[nosep]
  \item Theorems, lemmas, corollaries, and displayed equations are cited by their numbers in~\cite{DY25}.
  \item The paper indexes codewords by
        $c_1,\dots,c_4$; Lean uses \lean{row0 C}--\lean{row3 C} (i.e.\ rows
        $0$--$3$). We use the two conventions interchangeably.
  \item Inline identifiers such as
        \code{lambda C eps} are rendered from the source; module names appear in
        \texttt{typewriter} (e.g.\ \leanfile{Reduction.lean}); line ranges link to
        the corresponding source on GitHub.
\end{itemize}

\section{Formalization Design of Binary Codes}
\label{sec:encoding}

This section shows how the main theorems are formulated in Lean.
\Cref{sec:codes-lean} introduces the code model; \cref{sec:main-thms} gives
the Lean statements of Theorems~1--5 for comparison with the paper.

Lean~4~\cite{Lean4} is a dependently typed proof assistant: mathematical
objects and propositions are types, and a proof is a program whose type is the
theorem statement. Lean's trusted kernel checks these typing judgments.
The development lives in the \texttt{N4Code} namespace and builds on
mathlib~\cite{mathlib} for finite types (\texttt{Fin~$n$} for bounded
indices, \texttt{Finset} for finite sets), big sums written as
\code{Finset.sum ...}, parity (\lean{Even}/\lean{Odd}), and real arithmetic
($\mathbb{R}$).
We use $\mathbb{N}$ and $\mathbb{R}$ for the natural and real numbers;
\code{A -> B} for functions; $\forall$ and $\exists$ for quantifiers; and
\lean{Prop} for propositions.
An expression \code{t : T} means that \code{t} has type \code{T};
\lean{def} and \lean{abbrev} introduce definitions, while \lean{lemma} and
\lean{theorem} mark proved statements.

\subsection{Binary Codes in Lean}
\label{sec:codes-lean}

The paper \cite{DY25} studies optimal finite-length binary block codes of
size four for a binary symmetric channel (BSC) with crossover probability
$\eps$. A code $C \subseteq \{0,1\}^n$ of blocklength $n$
is optimal if it maximizes the correct-decoding probability
$\lambda_C(\eps)$ under maximum-likelihood decoding.

Mathematically, a binary $(n,4)$ code is a set of four distinct words
$c_1,\dots,c_4\in\{0,1\}^n$.
The paper indexes rows by $1,2,3,4$ and Lean by $0$--$3$; we use the two
conventions interchangeably.

\subsubsection{Basic Definitions of Codes}

The paper's proofs manipulate columns of the $4\times n$ matrix whose rows
are those codewords, so our formalization takes columns as the primary object:
\begin{Verbatim}[frame=single,breaklines=true,breakanywhere=true]
abbrev Word (n : ℕ) := Fin n → Bool
abbrev Column := Fin 4 → Bool
abbrev Code (n : ℕ) := Fin n → Column
\end{Verbatim}
Here \code{Fin k} is mathlib's type $\{0,1,\dots,k-1\}$;
a \lean{Column} is four bits indexed by \code{Fin 4};
and \code{Code n} is an $(n,4)$ code given by its $n$ columns.
The $j$-th codeword is \code{row C j : Word n}
($j=0,1,2,3$), also abbreviated \code{row0 C},\ldots,\code{row3 C}.

The type \code{Code n} is slightly more general than the notion of an
$(n,4)$ code used in~\cite{DY25}. In the paper, an $(n,4)$ code is a
four-element subset of $\{0,1\}^n$, and hence its four codewords are
pairwise distinct. In contrast, \code{Code n} represents an arbitrary
$4\times n$ binary matrix and therefore may also have repeated rows.
We deliberately work in this enlarged space because it simplifies the
intermediate column transformations used throughout the formalization.
Pairwise distinctness of the four rows is recorded separately in Lean by the
predicate \code{DistinctRows}.

Weight and distance on words can be defined as follows:
\begin{Verbatim}[frame=single,breaklines=true,breakanywhere=true]
def bitXor {n : ℕ} (x y : Word n) : Word n :=
  fun i => Bool.xor (x i) (y i)
def hammingWeight {n : ℕ} (x : Word n) : ℕ := 
    ∑ i ∈ Finset.univ, if x i = true then 1 else 0
def hammingDist {n : ℕ} (x y : Word n) : ℕ :=
  hammingWeight (bitXor x y)
\end{Verbatim}
Distance from a received word $y$ to codeword $c_j$, and to the code, are \lean{dRow} and \lean{dCode} respectively.
\begin{Verbatim}[frame=single,breaklines=true,breakanywhere=true]
def dRow {n : ℕ} (C : Code n) (j : Fin 4) (y : Word n) : ℕ :=
  hammingDist (row C j) y
def dCode {n : ℕ} (C : Code n) (y : Word n) : ℕ :=
  min (hammingDist (row0 C) y)
    (min (hammingDist (row1 C) y)
      (min (hammingDist (row2 C) y) (hammingDist (row3 C) y)))
\end{Verbatim}

\subsubsection{Performance of Codes}
On a BSC with crossover $\eps\in(0,1/2)$, the average correct-decoding
probability under ML decoding is
\begin{equation*}
  \lambda_C(\eps)
  \;:=\;
  \tfrac14\sum_{y\in\{0,1\}^n}
  (1-\eps)^{n-\dcode(y)}\,\eps^{\dcode(y)},
\end{equation*}
formalized as \code{lambda C eps} (\leanfile{Definitions.lean}, \leanat{Definitions.lean}{110}{111}).
The per-word contribution
$$w_{n,\eps}(d):=(1-\eps)^{n-d}\eps^d,$$ is formalized as
\lean{weight n eps d} (\leanfile{Performance.lean},
\leanat{Performance.lean}{23}{23}). With the counts $\alpha_C(d)$ of the number of $y$ with $d_C(y)=d$, we have
\begin{equation*}
\lambda_C(\eps)=\tfrac14\sum_d \alpha_C(d)
w_{n,\eps}(d).
\end{equation*}

Optimality and universal comparison are then
\begin{Verbatim}[frame=single,breaklines=true,breakanywhere=true]
def OptimalAt {n : ℕ} (C : Code n) (ε : ℝ) : Prop :=
  ∀ D : Code n, lambda C ε ≥ lambda D ε
def UniversalBetter {n : ℕ} (C₁ C₂ : Code n) : Prop :=
  ∀ ε : ℝ, 0 < ε → ε < 1 / 2 → lambda C₁ ε ≥ lambda C₂ ε
def UniversalStrictBetter {n : ℕ} (C₁ C₂ : Code n) : Prop :=
  ∀ ε : ℝ, 0 < ε → ε < 1 / 2 → lambda C₁ ε > lambda C₂ ε
def UniversalEqual {n : ℕ} (C₁ C₂ : Code n) : Prop :=
  ∀ ε : ℝ, 0 < ε → ε < 1 / 2 → lambda C₁ ε = lambda C₂ ε
\end{Verbatim}
so \code{OptimalAt C ε} is optimality at a fixed $\eps$, while
\code{UniversalBetter C' C} means $\lambda_{C'}(\eps)\ge\lambda_C(\eps)$
for every $\eps\in(0,1/2)$ (with a strict variant
\lean{UniversalStrictBetter} and an equality variant
\lean{UniversalEqual}).

Although \texttt{OptimalAt} quantifies over the enlarged type
\texttt{Code n}, this does not change the optimal-code problem for
$n\ge2$.  We prove that every code matrix with repeated rows is
universally strictly dominated by one with pairwise distinct rows.
Consequently, every \texttt{OptimalAt} code has pairwise distinct rows.
Thus the additional objects admitted by \texttt{Code n} can occur only
as intermediate nonoptimal matrices, and the optimal elements of the
Lean model coincide with genuine $(n,4)$ codes as defined in
\cite{DY25}.

\subsubsection{Equivalence and Code Classes}
Two codes are equivalent if one is obtained from the other by row
permutation, column permutation, and flipping all bits in some columns:
\begin{Verbatim}[frame=single,breaklines=true,breakanywhere=true]
def Equivalent {n : ℕ} (C C' : Code n) : Prop :=
  ∃ ρ : Equiv (Fin 4) (Fin 4), ∃ p : Equiv (Fin n) (Fin n),
    ∃ f : Fin n → Bool,
      ∀ t : Fin n,
        C' (p t) = rowPermute ρ
          (if f t then flipCol (C t) else C t)
\end{Verbatim}
Equivalence preserves $\lambda_C$, so classification statements are
naturally up to \lean{Equivalent}.

Each column has a type given by reading its four bits
as a binary integer in $\{0,\dots,15\}$. Note that the paper's column-type bit order is the
reverse of this natural row indexing: Bit $j$ of a column is row $(4-j)$'s entry, with row~$0$ most significant.
Named columns \lean{col0}, \lean{col1}, \lean{col3}, \lean{col5}, \lean{col6}, \lean{col7}
are the standard representatives of those types. The paper defines $|i|_C$ (or $|i|$) as the number of columns of type~$i$ in a code $C$:
\begin{Verbatim}[frame=single,breaklines=true,breakanywhere=true]
def colVal (c : Column) : ℕ :=
  ∑ j : Fin 4, if c j then 2 ^ (3 - j.val) else 0
def count {n : ℕ} (C : Code n) (i : ℕ) : ℕ :=
  ∑ t : Fin n, if colVal (C t) = i then 1 else 0
\end{Verbatim}

The sum of counts over a set $s$ of types is
\code{totalCounts C s}.
Linear codes and the three nonlinear classes (Class-I, Class-II, Class-III) used in the reductions are
predicates on these counts:
\begin{Verbatim}[frame=single,breaklines=true,breakanywhere=true]
def IsLinear {n : ℕ} (C : Code n) : Prop :=
  (∀ t : Fin n, colVal (C t) = 0 ∨ colVal (C t) = 3 ∨
    colVal (C t) = 5 ∨ colVal (C t) = 6) ∧
    ((count C 3 > 0 ∧ count C 5 > 0) ∨
      (count C 3 > 0 ∧ count C 6 > 0) ∨
      (count C 5 > 0 ∧ count C 6 > 0))
def ClassI {n : ℕ} (C : Code n) : Prop :=
  Odd (count C 1) ∧
    ((Even (count C 3) ∧ Even (count C 5) ∧ Even (count C 6)) ∨
      (Odd (count C 3) ∧ Odd (count C 5) ∧ Odd (count C 6))) ∧
    totalCounts C {1, 3, 5, 6} = n
def ClassII {n : ℕ} (C : Code n) : Prop :=
  count C 1 > 0 ∧ totalCounts C {1, 3, 5, 6} = n ∧
    ((Even (count C 1) ∧ Even (count C 3) ∧
        Odd (count C 5) ∧ Odd (count C 6)) ∨
      (Even (count C 1) ∧ Odd (count C 3) ∧
        Even (count C 5) ∧ Even (count C 6)))
def ClassIII {n : ℕ} (C : Code n) : Prop :=
  totalCounts C {1, 3, 5, 6, 7} = n ∧
    ((count C 1 = 1 ∧ count C 7 = 1 ∧ count C 6 = 0 ∧
        Even (count C 3) ∧ Odd (count C 5)) ∨
      (count C 1 = 1 ∧ count C 5 = 0 ∧ count C 7 = 0 ∧
        Odd (count C 3) ∧ Odd (count C 6)))
\end{Verbatim}
In short: \lean{IsLinear} allows only types $0,3,5,6$ with at
least two of $|3|,|5|,|6|$ positive;
\lean{ClassI}/\lean{ClassII} use only
$\{1,3,5,6\}$ under the stated parity patterns; and
\lean{ClassIII} additionally allows type~$7$ under one of two
count patterns.
The canonical linear representative with $|3|=n_3$, $|5|=n_5$, $|6|=n_6$
is \code{linearCode n3 n5 n6}.

\subsection{Statements of Main Theorems}
\label{sec:main-thms}

Using the vocabulary of \cref{sec:codes-lean}, the five
headline classification theorems of~\cite[\S2]{DY25} are stated as follows.

\subsubsection{Theorem~1 ({\normalfont\texttt{optimal\_in\_linear\_or\_class1}})}
For every $n\ge 2$ and every $\eps\in(0,1/2)$ there
exists a code that is optimal at $\eps$ and is linear or Class-I. 
The theorem is stated in \leanfile{Reduction.lean}, \leanat{Reduction.lean}{5686}{5688}:
\begin{Verbatim}[frame=single,breaklines=true,breakanywhere=true]
theorem optimal_in_linear_or_class1 (n : ℕ) (hn : 2 ≤ n) :
    ∀ ε : ℝ, 0 < ε → ε < 1 / 2 →
      ∃ C : Code n, (IsLinear C ∨ ClassI C) ∧ OptimalAt C ε
\end{Verbatim}
The hypothesis $n\ge 2$ makes the paper's standing blocklength assumption
explicit; the conclusion is per-$\eps$.

\subsubsection{Theorem~2 ({\normalfont\texttt{linear\_opt\_*}})}
Theorem~2 of the paper classifies optimal linear
$(n,4)$ codes by $n\bmod 3$ and states that each listed code is universally
strictly better than every other non-equivalent linear code.
The formalization gives four declarations in \leanfile{Linear.lean}, one per
residue class:
\begin{itemize}[nosep]
  \item \lean{linear_opt_residue2} (\leanat{Linear.lean}{9835}{9838}): $n=3k-1$,
        $C(k,k,k-1)$ is optimal;
  \item \lean{linear_opt_n3} (\leanat{Linear.lean}{10331}{10336}): $n=3$,
        $C_A$, $C(1,1,1)$, $C(1,2,0)$ are optimal;
  \item \lean{linear_opt_residue0} (\leanat{Linear.lean}{10625}{10633}): $n=3k$
        ($k\ge 2$), $C(k+1,k+1,k-2)$ and $C(k+1,k,k-1)$ are optimal;
  \item \lean{linear_opt_residue1} (\leanat{Linear.lean}{10715}{10723}): $n=3k+1$,
        $C(k+1,k,k)$ and $C(k+2,k,k-1)$ are optimal.
\end{itemize}
Each concludes \texttt{UniversalStrictBetter} (strict universal optimality).

The paper's closing domination sentence is imprecise when several codes are
listed for the same~$n$ ($n=3$; $n=3k$, $k\ge2$; $n=3k+1$): they are
pairwise non-equivalent yet can have equal~$\lambda$ (e.g.\ at $n=3$,
$\lambda(C(1,1,1))=\lambda(C(1,2,0))=\lambda(C_A)$), so no listed code
strictly dominates all the others.
The claim is correct only if ``not equivalent to them'' is read
collectively: a competitor~$D$ must be non-equivalent to \emph{every}
listed optimum before \emph{each} listed code strictly dominates~$D$.
The $n=3k-1$ case lists a single optimum $C(k,k,k-1)$ and needs no such
collective reading.

The $n=3$ declaration is representative; it requires a separate
non-equivalence hypothesis for each listed code:
\begin{Verbatim}[frame=single,breaklines=true,breakanywhere=true]
theorem linear_opt_n3 :
    ∀ D : Code 3, IsLinear D →
      ¬ Equivalent CA D → ¬ Equivalent (linearCode 1 1 1) D →
        ¬ Equivalent (linearCode 1 2 0) D →
        UniversalStrictBetter CA D ∧
          UniversalStrictBetter (linearCode 1 1 1) D ∧
            UniversalStrictBetter (linearCode 1 2 0) D
\end{Verbatim}
Here \lean{CA} is the code with columns $(3,5,0)$.

The same pattern appears in \lean{linear_opt_residue0} and
\lean{linear_opt_residue1} (two listed optima, two
non-equivalence hypotheses); \lean{linear_opt_residue2} states
the single-optimum case directly.
All four are proved in \leanfile{Linear.lean}.

\subsubsection{Theorem~3 ({\normalfont\texttt{universal\_strict\_better\_of\_not\_class}})}
In the paper, Theorem~3 is stated as follows:
\begin{quotation}
  For $(n,4)$ codes with $n>3$, the following
property universally holds: the set comprising the optimal linear codes,
Class-I codes, Class-II codes, and their equivalent codes contains all the
optimal codes.
\end{quotation}
\begin{sloppypar}
The initial Lean formalization used a per-$\eps$ classification instead:
\lean{optimal_equivalent_linear_class1_class2}
(\leanfile{Nbig.lean}, \leanat{Nbig.lean}{2206}{2210}):
\end{sloppypar}
\begin{Verbatim}[frame=single,breaklines=true,breakanywhere=true]
theorem optimal_equivalent_linear_class1_class2 (n : ℕ) (hn : n > 3) :
    ∀ ε : ℝ, 0 < ε → ε < 1 / 2 → ∀ C : Code n,
      OptimalAt C ε →
        ∃ C' : Code n, Equivalent C C' ∧
          (IsLinear C' ∨ ClassI C' ∨ ClassII C')
\end{Verbatim}
However, this formalization does not accurately state the meaning of
``universally'' in the theorem in the paper.

Comparing the paper and Lean proofs led us to a sharper universal
statement: any code outside the three classes has another code that is
strictly better at \emph{every} crossover probability.
We formalize this as
\lean{universal_strict_better_of_not_class}
(\leanfile{Nbig.lean}, \leanat{Nbig.lean}{2193}{2199}):
\begin{Verbatim}[frame=single,breaklines=true,breakanywhere=true]
theorem universal_strict_better_of_not_class (n : ℕ) (hn : n > 3) (C : Code n)
    (hnot : ∀ C' : Code n, Equivalent C C' →
      ¬ (IsLinear C' ∨ ClassI C' ∨ ClassII C')) :
    ∃ D : Code n, UniversalStrictBetter D C
\end{Verbatim}
The theorem \lean{optimal_equivalent_linear_class1_class2} is a corollary of this one.

\subsubsection{Theorem~4 ({\normalfont\texttt{universal\_strict\_better\_of\_not\_linear}})}
Theorem~4 gives a sufficient condition under which all optimal $(n,4)$
codes are equivalent to linear codes.
As with Theorem~3, we use a strengthened universal form: under the
Class-I replacement condition, any code not equivalent to a linear code
admits another that is universally strictly better.
The formal statement is
\lean{universal_strict_better_of_not_linear}
(\leanfile{Nbig.lean}, \leanat{Nbig.lean}{2218}{2224}):
\begin{Verbatim}[frame=single,breaklines=true,breakanywhere=true]
theorem universal_strict_better_of_not_linear (n : ℕ) (hn : n > 3)
    (hcond : ∀ C : Code n, ClassI C → count C 1 ≥ 3 →
      ∀ t : Fin n, C t = col1 →
        UniversalBetter (replaceColumn C t (argminType C)) C) :
    ∀ C : Code n,
      (∀ C' : Code n, Equivalent C C' → ¬ IsLinear C') →
        ∃ D : Code n, UniversalStrictBetter D C
\end{Verbatim}
Here, a one-column change is \code{replaceColumn C t s'};
\code{argminType C} is the type among $\{3,5,6\}$ of minimal
count (ties broken in that order).

When Theorem~4's hypothesis holds, Theorem~2 implies that optimal linear
codes are universally strictly better than every non-equivalent competitor.

\subsubsection{Theorem~5 ({\normalfont\texttt{optimal\_codes\_small\_n}})}

Theorem~5 lists all optimal $(n,4)$ codes for $2\le n\le 8$.
The strengthened Lean statement is \lean{optimal_codes_small_n}
(\leanfile{FiniteN.lean}, \leanat{FiniteN.lean}{2336}{2342}): for $n\neq3$, a code not equivalent to a
linear code has a universally strictly better code; for $n=3$, a code not
equivalent to one of the five codes has a universally strictly better code.

For $4\le n\le 8$, Theorem~5 follows from verifying Theorem~4's
hypothesis.
The Lean statement is \lean{n4to8_strict_better_of_not_linear}
(\leanfile{FiniteN.lean}, \leanat{FiniteN.lean}{2220}{2225}): if $C$ is not equivalent to a linear
code, then some code is universally strictly better than $C$.
\begin{Verbatim}[frame=single,breaklines=true,breakanywhere=true]
theorem n4to8_strict_better_of_not_linear (n : ℕ) (hn4 : n > 3) (hn8 : n ≤ 8) :
    ∀ C : Code n,
      (∀ C' : Code n, Equivalent C C' → ¬ IsLinear C') →
        ∃ D : Code n, UniversalStrictBetter D C
\end{Verbatim}

For $n=2$, Theorem~4 does not apply because it requires $n>3$.
The Lean statement is
\lean{n2_strict_better_of_not_linear}
(\leanfile{FiniteN.lean}, \leanat{FiniteN.lean}{260}{260}). 

For $n=3$, the Lean statement is
\lean{n3_strict_better_of_not_InOptimal3}
(\leanfile{FiniteN.lean}, \leanat{FiniteN.lean}{2253}{2255}): if a code is not equivalent to one of the
five codes in \lean{InOptimal3} (\lean{code135}, \lean{code136},
\lean{CA}, \code{linearCode 1 1 1},
\code{linearCode 1 2 0}), then some code is universally strictly
better than it. 
The converse is also formalized by
\lean{n3_representatives_optimal}: each of the five representatives
is universally optimal among all $(3,4)$ codes.

\section{Proof Formalization Development}
\label{sec:ai-assisted-formalization}

The initial code produced by the assistant consisted of definitions
together with a catalogue of statements to prove (each carrying a
\texttt{sorry}).
AI agents then worked in parallel to discharge those \texttt{sorry}s.
Development was tracked in a git repository, with each agent operating on
its own clone (worktree) of the codebase.
All \texttt{sorry}s were removed within approximately 48 hours, with only
light human intervention.
The resulting library contains no \texttt{sorry}, \texttt{admit}, or
\texttt{axiom} declarations; every headline theorem is proved, and
\texttt{lake build} succeeds.

A \texttt{sorry}-free build is not enough by itself: the Lean kernel checks
proofs relative to the axioms they depend on, so those axioms must be
audited as well.
\Cref{sec:audit} records how the development was made to depend only on
Lean's standard base axioms: the \texttt{native\_decide} trust axioms that
appeared in an earlier draft were eliminated by replacing them with
\texttt{decide}.

During the audit, we also found redundant proofs---for example, distances
between codewords were computed differently in different modules.
We guided the AI to rewrite them for consistency and brevity; the results are
summarized in \Cref{sec:simplifications} and do not change any theorem
statement.
The simplified formalization totals 37{,}561 lines across fourteen modules;
Table~\ref{tab:modules} gives the size and main content of each.

Finally, we extract the generic core of the formalization that is independent of this specific $(n,4)$ classification problem. We summarize the core in \Cref{sec:reusable-core}.

\subsection{Axiom Audit}
\label{sec:audit}
\subsubsection{Trusted Axioms}
\label{sec:axioms}

A successful build with no \texttt{sorry} is not enough: the proofs must
depend only on an explicit allowlist of axioms.
The module \texttt{N4Code/AxiomCheck.lean} runs \texttt{\#print axioms} on
every headline theorem, and \texttt{scripts/axioms\_check.sh} checks that
each axiom list is contained in the allowlist (any other axiom---in
particular \texttt{sorryAx}---fails the check).

The first class of permitted axioms is Lean's standard base axioms
\texttt{propext}, \texttt{Quot.sound}, and \texttt{Classical.choice}:
\begin{itemize}[nosep]
  \item \texttt{propext} (propositional extensionality): equivalent
        propositions are equal, i.e.\ $(P\leftrightarrow Q)\to P=Q$;
  \item \texttt{Quot.sound}: if $x\sim y$ under a relation used to form a
        quotient type, then the images of $x$ and $y$ in the quotient are
        equal;
  \item \texttt{Classical.choice}: every nonempty type has an element
        (Hilbert's $\varepsilon$-operator), which yields classical logic
        (including the law of excluded middle).
\end{itemize}
These are the default axioms of Lean's mathematics library; they are not
specific to this formalization.

\subsubsection[decide vs.\ native\_decide]{From \texttt{native\_decide} to \texttt{decide}}
\label{sec:native_decide_audit}
A proposition \(P\) is \emph{decidable} in Lean when an instance of
\texttt{Decidable P} supplies either a proof of \(P\) or a proof of
\(\neg P\). The function \texttt{decide} maps such a proposition to a
boolean, and the tactic \texttt{decide} closes a goal \(P\) by reducing
\texttt{decide P} to \texttt{true} inside the Lean kernel. Since the
reduction is pure kernel computation, this path adds no axiom beyond the
standard base axioms.

The tactic \texttt{native\_decide} takes a different path: it compiles the
same decision procedure to native code, executes it in the Lean runtime,
and---if the result is \texttt{true}---records a \emph{trust axiom}
asserting the theorem. In \texttt{\#print axioms} these trust facts appear
as per-lemma names of the form
\texttt{<lemma>.\_native.native\_decide.ax\_*}, for example
\path{two_zero_columns._native.native_decide.ax_1_10}. Accepting such
axioms makes the Lean compiler and runtime part of the trusted computing
base.

In the released version we eliminate this extra trust assumption entirely.
Every \texttt{native\_decide} in \texttt{N4Code} is replaced by
\texttt{decide}; the finite checks---concrete column-type arithmetic
(\texttt{colVal}, \texttt{testBit}), exhaustive searches over
$\{0,1\}^n$ for $n\le 3$, and equivalence witnesses on fixed matrices---are
small enough to be closed by kernel reduction. Where the default kernel
recursion limit was too small, the relevant module sets
\texttt{set\_option maxRecDepth 1000000}.
After the replacement, \texttt{N4Code/AxiomCheck.lean} reports only
\texttt{propext}, \texttt{Classical.choice}, and \texttt{Quot.sound} for
every headline theorem, and \texttt{scripts/axioms\_check.sh} no longer
permits \texttt{Lean.ofReduceBool} or
\texttt{*.\_native.native\_decide.ax\_*}.

The former distribution of \texttt{native\_decide} declarations, retained
only as development history, is shown in Table~\ref{tab:native-decl}.

\begin{table}[ht]
  \centering
  \caption{Former declarations carrying \texttt{native\_decide}, by module
  (all replaced by \texttt{decide} in the released version).}
  \label{tab:native-decl}
  \small
  \begin{tabular}{@{}lrlr@{}}
  \toprule
  \textbf{Module} & \textbf{Declarations} & \textbf{Module} & \textbf{Declarations} \\
  \midrule
  \texttt{Reduction} & 91 & \texttt{ZeroColumn} & 10 \\
  \texttt{FiniteN} & 33 & \texttt{ClassI} & 17 \\
  \texttt{Nbig} & 24 & \texttt{Definitions} & 7 \\
  \texttt{TwoColumn} & 10 & \texttt{Basic} & 2 \\
  \texttt{Linear} & 17 & \texttt{Performance} & 4 \\
  \bottomrule
  \end{tabular}
  \end{table}

In the earlier draft, every one of the 215
\texttt{native\_decide}-bearing declarations (lemmas, theorems, and the
relevant \texttt{example} blocks) carried a marker line immediately above
its docstring:
\begin{Verbatim}[frame=single,fontsize=\footnotesize,breaklines=true,breakanywhere=true]
-- native_decide: Contentful · n=3 · checked 2026-08-24
\end{Verbatim}
The marker records the audit group (\texttt{Mechanical} or
\texttt{Contentful}), the block length ($n=2$, $n=3$, or \texttt{any} when the statement is generic), and
the date of the check.
After the replacement these marker lines are retained but now read
\texttt{-- decide:}; the same declarations can be enumerated with
\begin{Verbatim}[frame=single,fontsize=\footnotesize,breaklines=true,breakanywhere=true]
rg -- "-- decide:" N4Code/
\end{Verbatim}
and checked, per marker, that the surrounding declaration is of the stated
group and block length.

\subsection{Simplifications}
\label{sec:simplifications}

These simplifications are internal proof cleanup; no
theorem statement or paper correspondence changed.

\begin{sloppypar}
\begin{itemize}[nosep]
  \item \textbf{Row-distance lemmas.}
        While re-examining the formalization, we replaced the repeated row-pair
        distance instances in \leanfile{Reduction.lean} and \leanfile{Nbig.lean} with one
        general lemma \lean{hammingDist_rows_of_types}: for a code whose column types
        lie in a finite set $S$, the Hamming distance between rows $i$ and $j$ is the
        weighted sum over $S$ of the type counts whose $i$-th and $j$-th bits differ.
        The specialized $\{1,3,5,6\}$, $\{1,3,5,7\}$, and $\{3,5,6\}$ statements are
        now short corollaries; the four full-type-set formulas
        \lean{hammingDist_row0_row2_eq}, \lean{hammingDist_row1_row2_eq},
        \lean{hammingDist_row2_row3_eq}, and \lean{hammingDist_row1_row3_eq} for
        \lean{Columns07} codes are one-line corollaries with $S =$ \lean{Finset.Icc 0 7}.
        The lemma, together with \lean{sum_indicator_of_types}, now lives in
        \leanfile{Basic.lean}; the redundant \lean{sum_colVal_indicator} helper in
        \leanfile{TwoColumn.lean} is removed.

  \item \textbf{Count-replacement lemmas.}
        The three general column-replacement lemmas \lean{count_replace_dec},
        \lean{count_replace_inc}, and \lean{count_replace_eq} now live in
        \leanfile{Compare.lean} next to \lean{sum_split_at}.  The 1-to-3 instances in
        \leanfile{Reduction.lean} and the zero-column instances in \leanfile{Linear.lean}
        and \leanfile{ZeroColumn.lean} are short corollaries of this family, instead of
        re-deriving the same sum-splitting argument.

  \item \textbf{Counts across involutive relabellings.}
        The eight \lean{count_swap12Code_*} and \lean{count_swap36Code_*} lemmas each
        re-derived the same per-type counting argument.  They are now one-line corollaries
        of \lean{count_involution_map} in \leanfile{Basic.lean} (counts transport across an
        involutive type relabelling), together with the \lean{swapVal12} / \lean{swapVal36}
        involutions and the \lean{colVal_swap12Code} / \lean{colVal_swap36Code} type maps.
        The same general fact now also absorbs \lean{count_swapRows01Code} and
        \lean{count_swapRows02Code} in \leanfile{Reduction.lean}.

  \item \textbf{Positive-count witnesses.}
        The three lemmas \lean{exists_col0_of_count_pos},
        \lean{exists_col1_of_count_pos}, and \lean{exists_col7_of_count_pos} were
        specialized copies of \lean{count_pos_iff_exists}; they are now one-liners
        combining that general fact with the matching \texttt{colVal\_eq\_*\_iff\_col*}
        lemma.

  \item \textbf{Class-I hypothesis reuse for the $3\leftrightarrow 6$ swap.}
        The col6 branches of \lean{class1_one} and \lean{class1_min} previously
        re-derived the Class-I hypotheses (\lean{htypes}, the parity conditions, and
        \lean{totalCounts}) on \lean{swap36Code C} by hand.  A new \lean{ClassI_swap36Code}
        lemma supplies the Class-I property of the swapped code, and \lean{classI_hyps}
        now unpacks it, so those branches work directly from \lean{ClassI C}.  The
        now-unused \lean{types_swap36Code} was removed.

  \item \textbf{Linear block counting.}
        The three \lean{linear_count_3/5/6} lemmas counted the blocks of a linear code,
        with \lean{linear_count_5} and \lean{linear_count_6} each re-proving the same
        offset bijection.
        We added the general \lean{card_fin_interval} and
        \lean{card_fin_suffix} helpers (next to the existing prefix helper
        \lean{card_fin_lt}), so the middle and suffix counts are now one-liners.

  \item \textbf{Unified \texttt{colVal} iff lemmas.}
        We consolidated the six \texttt{colVal\_eq\_*\_iff\_col*} lemmas into two general
        facts, \lean{colVal_colOfNat} and \lean{colVal_eq_iff_colOfNat}.  The
        named-column instances for types $0,1,3,5,6,7$ are now short corollaries.
\end{itemize}
\end{sloppypar}

\subsection{Reusable Generic Core}
\label{sec:reusable-core}

The simplification pass described above reworked the development around small,
reusable lemmas in addition to shortening proofs.  These components form two
tiers.

\begin{sloppypar}
\begin{itemize}[nosep]
  \item \emph{Generic infrastructure.}  The foundational modules
        \texttt{Definitions}, \texttt{Basic}, \texttt{Performance}, and
        \texttt{Comparison} (together with \texttt{Compare} for the one-column
        machinery) provide, without any Class-I/II/III dependence, the Hamming
        geometry and code model, the type-number bijection
        (\lean{colVal_colOfNat} and \lean{colVal_eq_iff_colOfNat}), equivalence
        invariance and transport (\lean{dRow_equiv}, \lean{dCode_equiv},
        \lean{alpha_equiv}, \lean{lambda_equiv},
        \lean{universalEqual_of_equivalent},
        \lean{universalStrictBetter_of_equivalent}), the BSC weight and order
        machinery, and the generic comparison engine (\lean{compare_bij}).
        The simplification also extracted the general counting and grouping
        facts \lean{sum_split_at}, \lean{count_replace_dec},
        \lean{count_replace_inc}, \lean{count_replace_eq},
        \lean{sum_indicator_of_types}, \lean{hammingDist_rows_of_types}, and
        \lean{count_involution_map}, together with the \lean{card_fin_lt},
        \lean{card_fin_interval}, and \lean{card_fin_suffix} block-counting
        helpers.
  \item \emph{Paper-specific but factored.}  Components tied to the 4-column
        Class-I setting---such as \lean{swapVal12}, \lean{swapVal36},
        \lean{colVal_swap12Code}, \lean{colVal_swap36Code},
        \lean{Equivalent_swap36Code}, \lean{ClassI_swap12Code},
        \lean{ClassI_swap36Code}, and \lean{classI_hyps}---are not generic, but
        they are now factored so the classification proofs reuse them instead of
        repeating ad-hoc case analyses.
\end{itemize}
\end{sloppypar}

The remaining modules build the classification on the generic tier and depend
only on it.  We expect the generic infrastructure---the distance-distribution
machinery, the cumulative-sum criterion (Corollary~21), and the
equivalence-invariance machinery---to be reusable by future formalizations in
information theory, complementing mathlib's existing coding-theory library.
The machine-checked paper-to-Lean dictionary and the consistency tooling keep
the correspondence auditable as the generic core is put to new uses.

\section{Notable Discrepancies}
\label{sec:discrepancies}

During verification we found that the formalization follows the paper
closely overall, but Lean code is not a direct translation of the prose:
formalization demands precise statements, and the paper's natural language is
not always precise, explicit, or correct.
The subsections below and Table~\ref{tab:discrepancies} catalogue the main
issues by theme; each subsection gives fuller detail.

\begin{table}[t]
  \centering
  \caption{Summary of discrepancies between the published paper and the
    formalization.}
  \label{tab:discrepancies}
  \small
  \begin{tabular}{@{}llp{0.48\textwidth}@{}}
  \toprule
  \textbf{Location} & \textbf{Type} & \textbf{Issue} \\
  \midrule
  \S II-A, type invariance & gap & Multiset invariance stated, not proved \\
  Cor.~10 & proof gap & Thm.~8 only for $w(c_3\oplus c_4)$ even \\
  Thm.~16, $n=3$ & strength & Vacuous for genuine codes; Lean adds \texttt{DistinctRows} \\
  Thm.~16, argmin & proof gap & $|5|,|6|$ minimizer cases incomplete \\
  Thm.~4 & proof gap & ``Equivalent to Class-II'' step \\
  Lem.~15 Case~1 & labeling & ``Class-II-b'' should be Class-III-b \\
  Thm.~3 proof & imprecision & Linear image: one of $|5|,|6|$ zero \\
  Thm.~16 proof, $n=2$ & infeasible-case error & The \(n=2\) case is impossible under the Class-I parity conditions\\
  Thms.~8, 11 & convention & Equality characterizations need types $0$--$7$ \\
  Thm.~12, eq.~(244) & convention & $\mathcal{Y}_3$ converse needs consistency \\
  Thm.~6(2) & convention & 0--15 column identification \\
  \bottomrule
\end{tabular}
\end{table}

\subsection[Code Type Invariance]{Code Types Under Equivalence}

The paper states that code types are invariant under equivalence: no
equivalent code should belong to two different classes among linear,
Class-I, and Class-II.
It offers no formal proof, but appeals to the following facts.
In Section~II-C, eqs.~(21)--(22) claim that for equivalent
codes $C$ and $C'$ with columns only in types $0$--$7$,
\[
  \{|1|,|2|,|4|,|7|\}_C = \{|1|,|2|,|4|,|7|\}_{C'}
  \quad\text{and}\quad
  \{|3|,|5|,|6|\}_C = \{|3|,|5|,|6|\}_{C'}
\]
as multisets. For example, flipping a type-1
column yields type~14, and swapping rows~1 and~4 turns that into type~7, so
a code with a type-1 column can be equivalent to one with a type-7 column, while preserving the multiset $\{|3|,|5|,|6|\}$.

Lean does not formalize the paper's multiset invariance directly; instead,
equivalence preserves the parity of column Hamming weights
(\lean{allEvenWeights_equiv}).  Since linear codes have only even-weight
columns (\lean{IsLinear_all_even}) while Class-I, Class-II, and Class-III
codes each have a type-1 column, the lemmas
\lean{no_cross_class_linear_class1}, \lean{no_cross_class_linear_class2},
and \lean{no_cross_class_linear_class3} in \leanfile{Performance.lean} show that
no linear code is equivalent to any of these three classes.  Separation among
Class-I, Class-II, and Class-III is not proved this way and would need a
finer invariant than column-weight parity alone.

\subsection{Gaps in the Column-Change Arguments}

Corollary~10 shows that a non-Class-I nonlinear code with
$|1|+|3|+|5|+|6|=n$ admits a linear or Class-I code with smaller~$|1|$.
The proof observes that at least one of
\[
  w(c_1 \oplus c_4) = |1|+|3|+|5|,\quad
  w(c_3 \oplus c_4) = |1|+|5|+|6|,\quad
  w(c_2 \oplus c_4) = |1|+|3|+|6|
\]
is even and invokes Theorem~8, claiming a better code with $|1|$ reduced
by one while preserving $|1|+|3|+|5|+|6|=n$.  Theorem~8, however, is proved
only when $w(c_3\oplus c_4)$ is even; extensions to the other row pairs appear
only in the remark after Theorem~8.

Lean closes this gap formally for the $w(c_1\oplus c_4)$- and
$w(c_2\oplus c_4)$-even cases that Theorem~8 leaves to its remark.
The latter is handled by \lean{one_bit_flip_1_to_5}; the former
uses \lean{swapRows01Code} (rows~$c_1\leftrightarrow c_2$, plus
column flips keeping types $0$--$7$) to reduce to the
$w(c_2\oplus c_4)$-even, $1\to5$ flip on an equivalent code while preserving
\lean{totalCounts} over $\{1,3,5,6\}$.
\lean{onepo_step} and
\lean{descent_to_linear_or_class1} in \leanfile{Reduction.lean}
assemble these into the full Corollary~10 formal proof.

\subsection{Proofs of Theorems~3 and~4}

The proofs of Theorems~3 and~4 contain gaps that Lean closes explicitly.

In the Class-III branch of Theorem~3, the paper says the linear image $C'$ can be chosen
``without type $\bspan{5}$ or $\bspan{6}$ columns''.
The intended reading is that exactly one of $|5|,|6|$ is zero: $|6|=0$ after
the Class-III-a map $(1,7)\to(3,5)$ and $|5|=0$ after the Class-III-b map
$1\to3$, with the other two linear types positive.
If both counts were zero, $C'$ would be all type~$3$.
Lean proves this precise form in \lean{classIII_not_optimal}.

In Lemma~15, Case~1, the three equality conditions of Theorem~8 are cited as
follows: (i)~$w(c_1\oplus c_3)$ and $w(c_2\oplus c_3)$ both odd (then $C$ is
Class-II); (ii)~$|1|=1$, $|5|=0$, $|3|$ and $|6|$ odd (where the paper says
``Class-II-b''); (iii)~$|1|=1$, $|6|=0$, $|3|$ and $|5|$ odd (where the paper
says ``equivalent to Class-II-b'').
Class-II-b requires $|1|$ even, so neither (ii) nor (iii) ($|1|=1$) is
Class-II-b; under Case~1 ($|2|=|4|=|7|=0$), (ii) is Class-III-b and (iii)
becomes Class-III-b after \lean{swapRows01Code} exchanges types
$5\leftrightarrow6$.
Lemma~15's conclusion is unaffected; only the two labels should read
``Class-III-b''.

Theorem~4 drives the $|1|\ge 3$ descent for Class-I codes: replace one
type-1 column by the argmin type $s\in\{3,5,6\}$ to obtain~$C'$ with
$|1|_{C'}=|1|-1$, then apply Lemma~14 to reach Class-I with $|1|$
reduced by~2.
The paper asserts that $C'$ ``is equivalent to a Class-II code'' before
invoking Lemma~14.
Lean closes this step in \lean{classI_descend_one}
(\leanfile{Nbig.lean}), with one branch per argmin type.
For $s=3$, $C'$ is Class-II directly.
For $s=5$, $C'$ is not Class-II, but \lean{swap12Code}
($c_2\leftrightarrow c_3$) is a code equivalence mapping $C'$ to Class-II.
For $s=6$, the needed $3\leftrightarrow6$ relabelling is likewise realizable
by \lean{Equivalent_swap36Code} (flip type-3/6 columns and swap
$c_1\leftrightarrow c_3$, as in \lean{swapRows02Code}).
In each branch the role-swapped code~$\tilde C$ is Class-II, Lemma~14
applies, and the $|1|\to|1|-2$ descent of Theorem~4 follows.

\subsection{Proof of Theorem~16}

Theorem~16 states that for a Class-I code with
$|1|=1$, replacing the type-1 column by the argmin type among
$\{3,5,6\}$ yields
$\lambda_{C'}>\lambda_C$ when $n\neq 3$ and $\lambda_{C'}=\lambda_C$ when
$n=3$.  

First, for a genuine $(n,4)$ code with pairwise distinct codewords, no
Class-I code with $|1|=1$ exists at $n=3$, so the $n=3$ equality clause is
vacuous.
Lean's \texttt{Code} type does not enforce distinct rows; for example, a
code with one type-$1$ column and two type-$3$ columns is Class-I with
$c_1=c_2$.
Here the argmin type is~$5$, so replacing the type-$1$ column yields the
linear code $C(2,1,0)$, which is strictly better.
The formalization therefore adds \texttt{DistinctRows} to
\lean{class1_one} (the Lean counterpart of Theorem~16).

Second, the proof of Theorem~16 states that for $2\le n\le 4$ and $n\ne 3$,
one has $|3|=|5|=|6|=1$. This statement should be restricted to $n=4$.
Indeed, since $C$ is Class-I with $|1|=1$, we have
$|3|+|5|+|6|=n-1$, and $|3|,|5|,|6|$ have the same parity. For $n=2$,
their sum would be $1$, which is impossible under the common-parity
condition; hence no such Class-I code exists. For $n=4$, their sum is $3$,
so the common-parity condition implies $|3|=|5|=|6|=1$. Thus the $n=2$
case is vacuous, and no separate $\mathcal{Y}_3^A$ argument is needed.
Accordingly, the previous discussion involving $|3|=0$,
$\{|5|,|6|\}=\{0,1\}$ and \texttt{classI\_count1\_not\_optimal} should be
removed.
The remark after Theorem~16 that code $C'$ is linear is likewise not part of
\lean{class1_one}, which states only the strict or equal $\lambda$
comparison.

The paper proves Theorem~16 when $|3|=\min\{|3|,|5|,|6|\}$ and, for the
remaining argmin types~5 and~6, asserts that the code can be transformed to
an equivalent one with $|3|=\min$, without giving the transformation.
Lean supplies these reductions in \leanfile{ClassI.lean}.
When $|5|=\min$, the row swap $c_2\leftrightarrow c_3$
(\lean{swap12Code}) is a code equivalence that exchanges
types~3 and~5 while fixing types~1 and~6, reducing to the $|3|=\min$ case
(\lean{class1_one_col5_strict}).
When $|6|=\min$, the same reduction is also realizable by equivalence---flip
every type-3/6 column and swap rows~$c_1$ and~$c_3$ (as in
\lean{swapRows02Code} in \leanfile{Reduction.lean}).


The $|3|=\min$ branch of Theorem~16 uses four lemmas in
\leanfile{ClassI.lean} that characterize $\mc Y_3$ and $\mc Y_5$ and give closed
forms for $\alpha^3_C(d)$ and $\alpha^5_C(d)$.
In the paper these statements divide by~$2$, and the quotients are integral only
under the Class-I parity condition $|3|\equiv|5|\equiv|6|\pmod 2$.
Lean restates them without divisions:
\begin{itemize}[nosep]
  \item \lean{Y5_iff_weights_count1}
        (\leanat{ClassI.lean}{455}{501}): paper eq.~(274), with the $|1|=1$
        simplifications (292)--(304).
        States $\mc Y_5$ membership by cleared weight inequalities only.

  \item \lean{Y3_iff_weights_count1}
        (\leanat{ClassI.lean}{521}{599}): paper eq.~(271), branches
        (281)--(287).
        Same $\mc Y_3$ split as the paper, with divisions removed and the
        boundary word assigned to the $\mc Y_3^A$ branch.

  \item \lean{alpha5_closed_count1}
        (\leanat{ClassI.lean}{619}{729}): paper eq.~(296).
        One triple sum over $(w_3,w_5,w_6)$ with indicator
        \lean{Y5WeightCond} and \lean{dRow3WeightEq}.

  \item \lean{alpha3_closed_count1}
        (\leanat{ClassI.lean}{751}{867}): paper eqs.~(283) and~(288).
        Merges the $\mc Y_3^A$ and $\mc Y_3^B$ sums into one triple sum with
        \lean{Y3WeightCondA} $\lor$ \lean{Y3WeightCondB} and
        \lean{dRow2WeightEq}.
\end{itemize}

\section{Final Remarks}
\label{sec:remarks}

We present a complete, kernel-checked Lean~4 formalization of the proofs from a research paper in coding theory. This note highlights the audit process, which is essential for ensuring the correctness of the formalization. Overall, with AI assistance, the formalization process can be completed without requiring significant human effort or expertise in formal mathematics. 

Why formalization? Three lessons stand out during our experience.
First, formalization gives a reproducible correctness check for a long,
complex mathematical argument.
Second, working at proof-script precision exposed imprecise hypotheses and
gaps in the human-written proofs, leading to strengthened statements and proofs.
Third, auditted formalization scripts can be used to guide the AI to generate more accurate and efficient proofs, and even new results.

One area that has not been explored in this note, but is interesting, is the use of AI to assist in extending and generalizing the current results on finite-length codes.  
This could be a promising direction for future research.

\printbibliography

\end{document}